\documentclass{article}
\usepackage{amsthm}
\usepackage[shortlabels]{enumitem}
\usepackage{amsmath}
\usepackage{amssymb}
\usepackage{ifthen}
\usepackage{color}
\usepackage{graphicx}
\usepackage{hyperref}         
\hypersetup{
  unicode=false,          
  pdftoolbar=true,        
  pdfmenubar=true,        
  pdffitwindow=false,     
  pdfstartview={FitH},    
  pdftitle={My title},    
  pdfauthor={Vikram Sharma},     
  pdfsubject={Subject},   
  pdfcreator={Vikram Sharma},   
  pdfproducer={Producer}, 
  pdfkeywords={keyword1} {key2} {key3}, 
  pdfnewwindow=true,      
  colorlinks=true,       
  linkcolor=red,          
  citecolor=green,        
  filecolor=magenta,      
  urlcolor=cyan           
}

         \newboolean{ABSTRACT}
        \setboolean{ABSTRACT}{true} 
        \setboolean{ABSTRACT}{false}
        \newcommand{\ifabstract}[2]{\ifthenelse{\boolean{ABSTRACT}}{#1}{#2}}
        \newcommand{\ifnabstract}[2]{\ifthenelse{\boolean{ABSTRACT}}{#2}{#1}}

\newcommand{\ignore}[1]{}        

\newtheorem{observ}{\sc Observation}
\newtheorem{lemma}{\sc Lemma}
\newtheorem{theorem}[lemma]{\sc Theorem}
\newtheorem{coro}[lemma]{\sc Corollary}
\newtheorem{propo}[lemma]{\sc Proposition}  

\newcommand{\blem}{\begin{lemma}}
  \newcommand{\elem}{\end{lemma}}
\newcommand{\bleml}[1]{\begin{lemma} \label{lem:#1}}
  \newcommand{\eleml}{\end{lemma}}
\newcommand{\blemT}[2]{\begin{lemma}[#1] \label{lem:#2}}
  \newcommand{\elemT}{\end{lemma}}

\newcommand{\bthm}{\begin{theorem}}
  \newcommand{\ethm}{\end{theorem}}
\newcommand{\bthml}[1]{\begin{theorem} \label{thm:#1}}
  \newcommand{\ethml}{\end{theorem}}
\newcommand{\bthmT}[2]{\begin{theorem}[#1] \label{thm:#2}}
  \newcommand{\ethmT}{\end{theorem}}

\newcommand{\bcor}{\begin{coro}}
  \newcommand{\ecor}{\end{coro}}
\newcommand{\bcorl}[1]{\begin{coro} \label{cor:#1}}
  \newcommand{\ecorl}{\end{coro}}
\newcommand{\bcorT}[2]{\begin{coro}[#1] \label{cor:#2}}
  \newcommand{\ecorT}{\end{coro}}

\newcommand{\bfac}{\begin{fact}}
  \newcommand{\efac}{\end{fact}}
\newcommand{\bfacl}[1]{\begin{fact} \label{fac:#1}}
  \newcommand{\efacl}{\end{fact}}
\newcommand{\bfacT}[2]{\begin{fact}[#1] \label{fac:#2}}
  \newcommand{\efacT}{\end{fact}}

\newcommand{\bobs}{\begin{observ}}
  \newcommand{\eobs}{\end{observ}}
\newcommand{\bobsl}[1]{\begin{observ} \label{obs:#1}}
  \newcommand{\eobsl}{\end{observ}}
\newcommand{\bobsT}[2]{\begin{observ}[#1] \label{obs:#2}}
  \newcommand{\eobsT}{\end{observ}}

\newcommand{\bpro}{\begin{propo}}
  \newcommand{\epro}{\end{propo}}
\newcommand{\bprol}[1]{\begin{propo} \label{pro:#1}}
  \newcommand{\eprol}{\end{propo}}
\newcommand{\bproT}[2]{\begin{propo}[#1] \label{pro:#2}}
  \newcommand{\eproT}{\end{propo}}

\newcommand{\bpf}{\begin{proof}}
  \newcommand{\epf}{\end{proof}}

\newcommand{\beq}{\begin{equation}}
\newcommand{\eeq}{\end{equation}}
\newcommand{\beql}[1]{\begin{equation}\label{eq:#1}}
\newcommand{\eeql}{\end{equation}}

\newcommand{\refeq}[1]{(\protect\ref{eq:#1})}

\newcommand{\CC}{{\mathbb C}}

\newcommand{\as}{\;\textcolor{blue}{\mathop{\mbox{\rm :=}}}\;}

\newcommand{\dd}{ ,\ldots , }

\newcommand{\ib}{\subseteq }

	\newcommand{\mmat}[2][ccccccccccccccccccccccccc]{\left[
		\begin{array}{#1}
		#2\\
		\end{array}\right]}

      \SetEnumerateShortLabel{r}{(\roman*)}
      \SetEnumerateShortLabel{P}{(P\arabic*)}
      \SetEnumerateShortLabel{R}{(R\arabic*)}

        \newcommand{\paren}[1]{\left( { #1 }\right)}        
        \newcommand{\abs}[1]{\left\lvert #1 \right\rvert}                

        \newenvironment{prog}{\begin{tabbing}
         xxx\=xxx\=xxx\=xxx\=xxx\=xxx\=xxx\=xxx\=xxx\=xxx\=xxx\=xxx\=xxx\=
        \kill\\}{
                \end{tabbing}}
        
        \newcommand{\bfp}{{\bf p}}

        \newcommand{\bfx}{{\bf x}}

        \newcommand{\sm}{\setminus}

        \newcommand{\efrac}[1]{\frac{1}{#1}}        

	\setboolean{ABSTRACT}{false}
	\setboolean{ABSTRACT}{true} 

	\newboolean{FIG_FLAG}
\setboolean{FIG_FLAG}{false}
\setboolean{FIG_FLAG}{true} 
	\newcommand{\figchoose}[2]{\ifthenelse{\boolean{FIG_FLAG}}{
	  #1}{#2}}
	\setboolean{FIG_FLAG}{true}	
	\setboolean{FIG_FLAG}{false}	

\newcommand{\sep}{\mathrm{sep}}
\newcommand{\comp}{\mathrm{\circ\,}}

\figchoose{

\newcommand{\sm}{\setminus}
}{}

\newcommand{\ii}{\jmath}
\newcommand{\ignoreannot}[2][\cored{ANNOTATION:}]{}%

\begin{document}
\title{Bounds on Distance to Variety in Terms of Coefficients of Bivariate Polynomials}
\date{}
\author{Vikram Sharma\\
The Institute of Mathematical Sciences, HBNI\\
Chennai, India 600113}
\maketitle
Let $f \in \CC[x]$ be a univariate polynomial of degree $d$ with roots $\alpha_1 \dd \alpha_d$. 
For a point $z\in \CC$, let $\sep(f,z) \as \min_i |z-\alpha_i|$.
Then for all points $z \in \CC$ we know that
the logarithmic derivative at $z$ is
        \beq
        {\frac{f'(z)}{f(z)}} = \sum_{i=1}^d \frac{1}{|z-\alpha_i|}.
        \eeq
and more generally for any $k \ge 1$ we have
        \beq
        {\frac{f^{(k)}(z)}{f(z)}} = \sum_{1 \le i_1 < i_2 < \cdots < i_k \le d} \frac{1}{(z-\alpha_{i_1}) (z-\alpha_{i_2}) \ldots (z-\alpha_{i_k})}.
        \eeq
Taking absolute value on both sides, applying triangular inequality on the RHS, and observing that
the number of terms on the RHS is ${d \choose k}$ and each is smaller than $1/\sep(f,z)$ we get
the following bound: for $k \ge 1$
        \beql{kder}
        \abs{\frac{f^{(k)}(z)}{f(z)}}^{1/k} \le \frac{d}{\sep(f,z)}.
        \eeql
Another way to interpret this bound is to state it as follows:
        \beql{sepbound}
        \sep(f, z) \le d \min_{1 \le k\le d} \abs{\frac{f(z)}{f^{(k)}(z)}}.
        \eeql
Similar bounds are also derived in \cite[p.~452]{henrici:applied-bk}.
In this short note, we will generalize this result to bivariate polynomial $f(x, y) \in \CC[x,y]$.
The analogue result will have the following form: the left hand side will be the distance of a
point $\bfp$ from the variety of $f$, and the RHS will consist of the total degree of $f$ and
a quantity dependent on the absolute values of $f$ and its partial derivatives evaluated at $\bfp$.
We first establish some notation. For $k \ge 0$, define
        \beql{pder}
        f_{i,k}(\bfp) \as \frac{\partial^k f(\bfp)}{\partial^i x \partial^{k-i}y}.
        \eeql
Let $D$ be the total degree of $f$, $V(f) \ib \CC^2$ be the variety of $f$, and
        \beql{sep}
        \sep(\bfp, V(f))  \as \inf_{\bfx \in V(f)} \|\bfp-\bfx\|
        \eeql
be the distance function to $V(f)$.

The idea for deriving the bound is as follows. Consider a point $\bfp=(\bfp_x, \bfp_y) \in \CC^2 \sm V(f)$.
In order to derive an upper bound on $\sep(\bfp, V(f))$, we will consider all the lines through $\bfp$.
These lines intersect the curve $f(x,y)=0$ at finitely many points that can be obtained as roots
of a univariate polynomial. For instance, consider the intersection of the line $x=\bfp_x$
with the curve $f=0$. Apply the upper bound in \refeq{sepbound} to the resulting univariate polynomial
we obtain that
        $$\sep(\bfp, V(f)) \le D \min_{1\le k \le D} \left|\frac{f(\bfp)}{f_{0,k}(\bfp)} \right|^{1/k}.$$
Similarly, considering the intersection of the line $y=\bfp_y$ with the curve $f=0$ we also get that
        $$\sep(\bfp, V(f)) \le \min_{1\le k \le D} \left|k! {D \choose k} \frac{f(\bfp)}{f_{k,0}(\bfp)} \right|^{1/k}.$$
How do we get the terms corresponding to the mixed partial derivatives? We consider all the lines
with slope $\tan \theta$, as $\theta$ varies from $0$ to $2\pi$, and take the minimum of the absolute
value of the corresponding roots over all $\theta$. Since this function is periodic in $\theta$, it makes
sense to use some tools from Fourier analysis. The remaining section develops this idea into full detail. 


Considering $f$ as a polynomial in $x$ with coefficients in $\CC[y]$,
from the local parameterization of algebraic curves \cite{walker:curves:bk},
we know that in a certain neighborhood of a point
$(x,y) \in \CC^2 \sm V(f)$ we can express
        \beql{fx}
        f(x,y) \as K \prod_{i=1}^{d(y)}(x - \alpha_i(y)),
        \eeql
where $\alpha_i$'s are holomorphic functions of $y$, the degree
$d(y) \leq \deg(f,x)$ depends on the $y$-coordinate, and $K \in \CC$
is some constant.
Differentiating both sides with respect to $x$ and factoring $f(x,y)$ from the RHS we obtain that
        \beql{dfx}
        f_{1,0}(x,y) =  f(x,y)\sum_{i=1}^{d(y)} \frac{1}{x-\alpha_i(y)},
        \eeql
and in general
        $$f_{k,0}(x,y) = f(x,y)\sum_{1 \le i_1 < i_2 < \cdots <i_k \le d(y)} \frac{1}{(x-\alpha_{i_1}(y)) \ldots (x-\alpha_{i_k}(y))}.$$
Following the argument used to derive \refeq{kder} in the univariate setting, we obtain that for any point $\bfp \in \CC^2$
        \beql{b1}
        \left| \frac{f_{k,0}(\bfp)}{f(\bfp)} \right|^{1/k} \le \frac{D}{\sep(\bfp,V(f))},
        \eeql
Note that if there is an asymptote at $y$ then $d(y) < \deg(f,x)$; also, if $d(y)=0$ then the bound
above on the partial derivatives trivially holds since all the partial derivatives vanish.

We want to derive a similar bound for the mixed partial derivatives $f_{i, k-i}(\bfp)$.
To obtain this, we change the coordinate system
and then consider the intersection with either the horizontal or vertical axis.
Consider the following change of coordinates:
\ignore{
        $$\vectwo{x}{y} \as \matrixtwo{\frac{e^{\ii
	\theta}}{2}}{\frac{\sqrt{3}}{2}}{\frac{\sqrt{3}}{2}}
		{\frac{-e^{-\ii \theta}}{2}} \vectwo{X}{Y}.$$
where $\ii=\sqrt{-1}$.
}
        $$\mmat{x\\y} \as 
		\efrac{\sqrt{2}} \mmat{ e^{\ii\theta} & e^{-\ii\psi}\\
					e^{-\ii\psi} &  -e^{-\ii\theta}}
		\mmat{ X \\ Y}
	 = 
		U \cdot \mmat{ X \\ Y}$$
where $\ii=\sqrt{-1}$ and $\theta,\psi$ are any angles; we will later set $\psi=0$.
Note that the matrix $U$ is unitary since
        $$UU^{\dagger} =
			\efrac{\sqrt{2}} \mmat{ e^{\ii\theta} & e^{-\ii\psi}\\
					e^{\ii\psi} &  -e^{-\ii\theta}}\cdot
			\efrac{\sqrt{2}} \mmat{ e^{-\ii\theta} & e^{\ii\psi}\\
					e^{-\ii\psi} &  -e^{\ii\theta}}
		= \mmat{1 & 0 \\ 0 & 1}.$$
Define 
	$$F(X,Y) \as f (U(X,Y)) = f\paren{\frac{e^{\ii \theta}X+ e^{-\ii\psi}Y}{\sqrt{2}},
		\frac{e^{\ii\psi}X-e^{-\ii \theta}Y}{\sqrt{2}}}.$$
By repeated applications of the chain rule of partial differentiation we know that
        $$
        F_{k,0}(X,Y) = \sum_{i=0}^k {k \choose i} (f_{i,k-i}
		\comp U(X,Y))  \paren{\frac{\partial x}{\partial X}}^i 
                           \paren{\frac{\partial y}{\partial X}}^{k-i}.
        $$
Observe that	
	$$\partial x/\partial X  = e^{\ii \theta}/\sqrt{2},\qquad
		\partial y/\partial X  = e^{\ii\psi}/\sqrt{2}.$$
Thus,
        \beql{fxk}
        F_{k,0}(X,Y)  =\sum_{i=0}^k {{k \choose i}} 
		(f_{i,{k-i}} \comp U(X,Y)) e^{\ii i\theta} 3^{\ii (k-i)\psi} {2^{-k}}.
        \eeql
Since total degree of $F$ is the same as the total degree of $f$,
it follows from \refeq{b1} that for a point $\bfp \in \CC^2$
        \beql{b2}
        \abs{\frac{F_{k,0}(U^{-1}(\bfp))}{F(U^{-1}(\bfp))}}^{1/k} 
                \le \frac{D}{\sep(U^{-1}(\bfp), V(F))}
                =\frac{D}{\sep(U^{-1}(\bfp), U^{-1}(V(f)))}
                =\frac{D}{\sep(\bfp,V(f))},
         \eeql
where the last step follows from the fact that  $U$ is a unitary transformation. 

Moreover, as $F(U^{-1}(\bfp)) = f(\bfp)$, 
from \refeq{fxk} and \refeq{b2} we obtain that for all 
choices of $\theta \in [-\pi, \pi]$
        \beql{ineq1}
        \abs{\sum_{i=0}^k {{k \choose i}} \frac{ f_{i,{k-i}}(\bfp)}{f(\bfp)}
		e^{\ii i\theta} e^{\ii (k-i)\psi} {2^{-k}}}^{1/k}
                \le \frac{D}{\sep(\bfp,V(f))}.
        \eeql
Let $P(\theta)$ be the function inside the absolute value on the LHS above.
Since it is a Fourier series in $\theta$, from Parseval's theorem we know that
        $$\sum_{i=0}^k \paren{{k \choose i}  \abs{\frac{f_{i,{k-i}}(\bfp)}{f(\bfp)}}
		e^{\ii (k-i)\psi} {2^{-k}}}^2
                = \frac{1}{2\pi} \int_{-\pi}^\pi |P(\theta)|^2 d\theta.$$
Substituting the upper bound \refeq{ineq1} on
$|P(\theta)|$ in the integral on the RHS above we further obtain that
        \beql{b4}
        \sum_{i=0}^k \paren{{k \choose i} 
		\abs{\frac{f_{i,{k-i}}(\bfp)}{f(\bfp)}} e^{\ii (k-i)\psi} {2^{-k}}}^2
                \le \paren{\frac{D}{\sep(\bfp,V(f))}}^{2k}.
        \eeql
Choosing $\psi=0$, we obtain
        $$
        \sum_{i=0}^k \paren{{k \choose i} 
		\abs{\frac{f_{i,{k-i}}(\bfp)}{f(\bfp)}} e^{\ii(k-i)\psi} {2^{-k}}}^2
                        > 
               \max_{j=0\dd k}\paren{
	       \abs{\frac{f_{j,{k-j}}(\bfp)}{f(\bfp)}}{2^{-k}}}^2
	 $$
Combining this with \refeq{b4} we obtain the following:

\bthml{dbound}
For a point $\bfp \in \CC^2$ that is not a zero of $f$
        \beql{pdiff}
        \max_{i=0 \dd k}\abs{\frac{f_{i,{k-i}}(\bfp)}{f(\bfp)}}^{1/k}
		<\frac{2D}{\sep(\bfp,V(f))}.
        \eeql
\ethml

The bound above can also be interpreted as an upper bound on the separation of a point $\bfp$
from the variety $V(f)$ in terms of the coefficients of the polynomial. Can a converse bound be given, i.e., 
a {\em lower bound} on the separation in terms of the coefficients. We next derive such a bound.

Suppose $f(0) \neq 0$, then we want to  derive a lower bound on $\sep(0, V(f))$ in terms of the coefficients.
Clearly, any $\bfx=(x,y)$ for which
        \beql{a00}
        |a_{0,0}| > \sum_{k \ge 1} \sum_{i=0}^k |a_{i,k-i}| |x|^i |y|^{k-i}
        \eeql
cannot be on the variety of $f$. Define
        $$\gamma \as \max_{1 \le k \le D} 
                \max_{0 \le i \le k} \paren{\frac{k!}{\binom{k}{i}} \left|\frac{a_{i,k-i}}{a_0}\right|}^{1/k},$$
where $D$ is the total degree. Then it follows that \refeq{a00} is equivalent to
        $$1 >  \sum_{k \ge 1} \frac{\gamma}{k!}\sum_{i=0}^k {k \choose i} |x|^i |y|^{k-i}
            =  \sum_{k \ge 1} \frac{\gamma}{k!} (|x|+ |y|)^k
            > \exp(\gamma \|\bfx\|_1) -1.$$
Therefore, if $\bfx$ is such that
$\|\bfx\|_1 \gamma < \ln 2$ then $|f(\bfx)| > 0$. In general, for any point $\bfp \in \CC^2$ we can 
apply the argument above to the shifted polynomial to obtain the following: if
        \beql{gammax}
        \gamma_f(\bfp) \as \max_{1 \le k \le D} 
                \max_{0 \le i \le k} { \left|\frac{f_{i,k-i}(\bfp)}{f(\bfp)}\right|}^{1/k},
         \eeql
then
        \beql{lsep}
        \sep(\bfp, V(f)) \ge \frac{\ln 2}{\sqrt{2}\gamma(\bfp)} \ge \frac{1}{3\gamma(\bfp)}.
        \eeql

Besides their intrinsic interest, such bounds are useful in analyzing the complexity of certain algorithms.
For instance, the bound given in \refeq{kder} has been useful in bounding the running time of certain
root isolation algorithms using the continuous amortization framework \cite{burr:contAmortization:16,sharma-batra:issac:2015}.
We expect the generalization given above to be useful in deriving similar bounds on the running
time of generalizations of corresponding algorithms that generally use subdivision (e.g., \cite{plantinga-vegter:isotopic:04}).

{\bf Acknowledgement:} The author is grateful to Chee Yap and Bernard Mourrain for their feedback
on earlier drafts of the results presented here.

\end{document}